\documentclass[12pt]{article}

\usepackage[slantedGreek]{mathptmx}
\usepackage{mathptmx}
\usepackage{float}
\usepackage[scaled=0.90]{helvet}
\usepackage{courier}
\usepackage{abbrevs}
\normalfont

\usepackage{hyperref}

\usepackage{graphicx}
\graphicspath{{illustrations/}}

\usepackage{amsmath}
\usepackage{amsthm,amssymb}
\usepackage{amsfonts}
\usepackage{graphics}
\usepackage{graphicx}
\usepackage[all]{xy}
\usepackage[usenames]{color}

\usepackage{natbib}
\bibpunct{(}{)}{,}{a}{,}{,}

\numberwithin{equation}{section}

\newtheorem{theo}{Theorem}[section]
\newtheorem{lem}[theo]{Lemma}
\newtheorem{cor}[theo]{Corollary}
\renewcommand{\thefootnote}{\fnsymbol{footnote}}

\begin{document}
\vspace{-0.5cm}
\begin{center}

{\huge\bf A nonparametric test for Cox processes\\
\vspace{0.2cm}}

\vspace{0.5cm}
Beno\^{\i}t \textsc{Cadre}

IRMAR, ENS Rennes, CNRS\\
Campus de Ker Lann\\
Avenue Robert Schuman,
35170 Bruz, France\\
\smallskip
{benoit.cadre@ens-rennes.fr}

\vspace{0.5cm}
Gaspar \textsc{Massiot} and Lionel \textsc{Truquet}

IRMAR, Ensai, CNRS\\
Campus de Ker Lann\\
Avenue Robert Schuman,
35170 Bruz, France\\
\smallskip
{massiot,truquet@ensai.fr}

\end{center}

\medskip

\begin{abstract} 
\noindent In a functional setting, we propose two test statistics to highlight the Poisson nature of a Cox process when $n$ copies of the process are available. Our approach involves a comparison of the empirical mean and the empirical variance of the functional data and can be seen as an extended version of a classical overdispersion test for counting data.
The limiting distributions of our statistics are derived using a functional central limit theorem for c\`adl\`ag martingales. We also study the asymptotic power of our tests under some local alternatives.  
Our procedure is easily implementable and does not require any knowledge of covariates.
 A numerical study reveals the good performances of the method. We also present two applications of our tests to real data sets.  \\

\end{abstract}

\renewcommand{\thefootnote}{\arabic{footnote}}

\setcounter{footnote}{0}

\section{Introduction}

Count process formulation is commonly used to describe and analyze many kind of data in sciences and engineering. A special class of such processes that researchers across in different fields frequently encounter are the so-called Cox processes or doubly stochastic Poisson processes. Compared to the standard Poisson process, the key feature of a Cox process is that its arrival rate is stochastic, depending on some covariate. In other words, if we let $T>0$ the observation period, $N=(N_t)_{t\in [0,T]}$ the Cox process and $\Lambda=(\Lambda(t))_{t\in[0,T]}$ the (stochastic) cumulative arrival rate then, conditioning on $\Lambda$, the distribution of $N$ is that of a Poisson process with cumulative intensity $\Lambda$. The benefit of randomness in the cumulative intensity lies in the fact that the statistician can take into account auxiliary informations, thus leading to a better model. For general references, we refer the reader to the monographies by Cox and Isham (1980), Karr (1991) and Kingman (1993). 

In actuarial sciences and risk theory for instance, the number of claims in the risk model may be represented by a Cox process. In this area, the central quantity is the ruin probability, that is the probability that the surplus of the insurer is below zero at some time (see e.g., Bj\"{o}rk and Grandell, 1988; Grandell, 1991; Schmidili, 1996). Cox process also appears in biophysics and physical chemistry (see e.g., Kou et al., 2005; Kou, 2008; Zhang and Kou, 2010). In these fields, experimental data consist of photon arrival times with the arrival rate depending on the stochastic dynamics of the system under study (for example, the active and inactive states of an enzyme can have different photon emission intensities); by analyzing the photon arrival data, one aims to learn the system's biological properties. Cox process data arise in neuroscience, to analyse the form of neural spike trains, defined as a chain of action potentials emitted by a single neuron over a period of time (see e.g., Gerstner and Kistler, 2002; Reynaud-Bourret et al., 2014). Finally mention astrophysics as another area where Cox process data often occur (see e.g., Scargle, 1998; Carroll and Ostlie, 2007). 

In general, it is tempting to associate to a model numerous covariates, and this possibly abusively. With this kind of abuse, one may consider a Cox process model though a Poisson process model is satisfactory. In this paper, we elaborate a nonparametric test statistic to highlight the Poisson nature of a Cox process. More precisely, based on i.i.d. copies of $N$, we construct a nonparametric test statistic for ${\bf H_0}$: $N$ is a Poisson process vs ${\bf H_1}$: $N$ is not a Poisson process. This setting of i.i.d. copies of $N$ is justified by the fact that in many situations, the  duration of observation is limited but the number of observed paths is large. 

Among the various possibilities to elaborate a test statistic devoted to this problem, one could estimate both functions $t\mapsto \mathbb E [N_t|\Lambda]$ and $t\mapsto \mathbb EN_t$ and test whether these functions are equal. However, this approach suffers from two main drawbacks, that is curse of dimensionality (whenever $\Lambda$ takes values in a high-dimensional space) and knowledge a priori on $\Lambda$. Another approach is to test whether time-jumps of $N$ are Poisson time-jumps; in this direction, we refer the reader to the paper by Reynaud-Bourret et al. (2014), in which a modified Kolmogorov-Smirnov statistic is used. 

In this paper, we elaborate and study a test statistic based on the observation that a Cox process is a Poisson process if, and only if its mean and variance function are equal. As we shall see, this approach leads to a very simple and easily implementable test.
  
The paper is organized as follows. In Section~\ref{section:test}, we first present the test statistic, then we establish asymptotic performances dedicated to the problem of ${\bf H_0}$ vs ${\bf H_1}$. The case of a local alternative is also considered. Section~\ref{section:simu} is devoted to a simulation study. An application to real data is presented in Section~\ref{section:data}. The proofs of our results are postponed to the three last sections of the paper.

\section{Tests for Cox processes}\label{section:test}

\subsection {Principle of the test}Ê
Throughout the paper, $T>0$ is the (deterministic) duration of observation, and $N=(N_t)_{t\in [0,T]}$ is a Cox process with aggregate process $\Lambda=(\Lambda(t))_{t\in [0,T]}$, such that $\mathbb E N_T^4<\infty$ and $\mathbb E N_t\ne 0$ for some $t\in ]0,T[$. Note that we do not assume here that $N$ has an intensity. 

We let $m$ and $\sigma^2$ the mean and variance functions of $N$, i.e. for all $t\in [0,T]$ : 
$$m(t)=\mathbb E N_t \mbox{ and } \sigma^2(t)={\rm var}(N_t).$$
Recall that for all $t\in [0,T]$ (see p. 66 in the book by Kingman, 1993) : 
\begin{equation}\label{king}
\sigma^2(t)=m(t)+{\rm var}(\mathbb E[N_t|\Lambda])=m(t)+{\rm var} \big(\Lambda(t)\big).
\end{equation} 
Hence, $\sigma^2(t)\ge m(t)$ that is, each $N_t$ is overdispersed. Moreover, if $m=\sigma^2$, then $\mathbb E[N_t|\Lambda]=\mathbb E N_t$ for all $t\in [0,T]$, thus $N$ is a Poisson process. As a consequence, $N$ is a Poisson process if, and only if $m=\sigma^2$. This observation is the key feature for the construction our test statistic, insofar the problem can be written as follows: 
$${\bf H_0}\, :\, \sigma^2=m \mbox{ vs } {\bf H_1} \, :\, \exists t\le T \mbox{ with } \sigma^2(t)>m(t).$$

From now on, we let the data $N^{(1)},\cdots,N^{(n)}$ to be independent copies of $N$. By above, natural test statistics are based on the process $\hat \sigma^2-\hat m=(\hat \sigma^2(t)-\hat m(t))_{t\in [0,T]}$, where $\hat m$ and $\hat \sigma^2$ are the empirical counterparts of $m$ and $\sigma^2$: 
$$\hat m(t)=\frac{1}{n} \sum_{i=1}^n N_t^{(i)} \mbox{ and } \hat \sigma^2(t)=\frac{1}{n-1}\sum_{i=1}^n \big(N_t^{(i)}-\hat m(t)\big)^2.$$

 In this paper, convergence in distribution of stochastic processes is intended with respect to the Skorokhod topology (see Chapter VI in the book by Jacod and Shiryaev, 2003). 

\begin{theo} \label{lim1} Let $B=(B_t)_{t\in\mathbb R_+}$ be a standard Brownian Motion on the real line. Under ${\bf H_0}$, $\hat \sigma^2-\hat m$ is a martingale and 
$$\sqrt{n}\, \big(\hat \sigma^2-\hat m\big)\stackrel{({\rm law})}\longrightarrow \big(B_{2m(t)^2}\big)_{t\le T}.$$
\end{theo}

As far as we know, the martingale property for $\hat \sigma^2-\hat m$ has not been observed yet. This property, which is interesting by itself, plays a crucial role in the derivation of the asymptotic result. 

\subsection{Testing ${\bf H_0}$ vs ${\bf H_1}$} 
Among the various possibilities of test statistics induced by $\hat \sigma^2-\hat m$, we shall concentrate in this paper on
$$\hat S_1 = \sup_{t\le T} \big(\hat \sigma^2(t)-\hat m(t)\big), \mbox{ and }
\hat S_2 = \int_0^T \big(\hat \sigma^2(t)-\hat m(t)\big) {\rm d} t.$$
These test statistics induced by the supremum and the integral are natural to test if a nonnegative function is equal to 0. Thus, they are chosen so as to respect the unilateral nature of the problem due to the fact that the alternative hypothesis may be written ${\bf H_1}$ : $\sigma^2(t)>m(t)$ for some $t\le T$. 

We now present the asymptotic properties of $\hat S_1$ and $\hat S_2$.

\begin{cor} \label{limite} Let ${\hat I}^{\, 2}=\int_0^T (T-t)\hat m(t)^2 {\rm d}t.$  

(i) Under ${\bf H_0}$, $$\sqrt{n}\frac{\hat S_1}{\hat m(T)} \, \stackrel{\rm{(law)}}\longrightarrow  \, |\mathcal N(0,2)|, \mbox{ and } \sqrt{n}\frac{\hat S_2}{\hat I} \, \stackrel{\rm{(law)}}\longrightarrow  \, \mathcal N(0,4).$$

(ii) Under ${\bf H_1}$, 
$$\sqrt{n}\frac{\hat S_1}{\hat m(T)}\, \stackrel{\rm prob.}\longrightarrow \, +\infty, \mbox{ and } \sqrt{n}\frac{\hat S_2}{\hat I}\, \stackrel{\rm prob.}\longrightarrow \, +\infty$$
\end{cor}
Hence, the test statistics $\hat S_1/\hat m(T)$ and $\hat S_2/\hat I$ define asymptotic tests with maximal power and the rejection regions for  tests of level $\alpha\in ]0,1[$ are: 
 \begin{equation}\label{RR} R_1(\alpha)=\Big\{\frac{\hat S_1}{\hat m(T)}\ge \sqrt{\frac{2}{n}}  \, q_{1-\alpha/2}\Big\} \mbox{ and } R_2(\alpha)=\Big\{\frac{\hat S_2}{\hat I}\ge \frac{2}{\sqrt{n}} \, q_{1-\alpha}\Big\},\end{equation}
 where for each $\beta\in ]0,1[$, $q_\beta$ is the $\mathcal N(0,1)$-quantile of order $\beta$. \\
 
\noindent {\bf Remark.} A close inspection of the proof of Theorem \ref{limite} reveals that a more general setting may be reached. Indeed, for the test of ${\bf H_0}$ vs ${\bf H_1}$, we only need to  assume that $N$ is in some class of overdispersed counting processes (i.e. ${\rm var} (N_t)\ge \mathbb E N_t$ for all $t\in [0,T]$) which satisfies the property : ${\rm var} (N_t)=\mathbb E N_t$ for all $t\in [0,T]$ if, and only if $N$ is a Poisson process. The archetype of such a class of counting processes is given by the Cox process, but it is also satisfied by other classes, such as some subclasses of Hawkes process for instance. In this direction, our test is more or less a functional version of the classical overdispersion test, that is used for testing the Poisson distribution of a sequence of count data (see for instance Rao and Chakravarti, 1956 or Bohning, 1994). Recall that overdispersion tests are widely used in actuarial science in the study of claims counts (e.g. Denuit et al., 2007) .

 \subsection{Local alternative}Ê
 By above, both tests defined by the rejection regions in \eqref{RR} are consistent. The aim of this section is to introduce a local alternative in view of a comparison of the asymptotic power of the tests. We refer the reader to Chapter 14 in the book by van der Vaart (1998) for a general overview on local alternatives.

 In this section, we assume in addition that the Cox process $N$ has an intensity $\lambda=(\lambda(t))_{t\in [0,T]}$, i.e. with probability 1, $\Lambda$ is absolutely continuous, and $$\Lambda(t)=\int_0^t \lambda(s){\rm d}s,\ \forall t\in [0,T].$$

We introduce a vanishing sequence of positive numbers $(d_n)_n$, and we consider the case of a local alternative defined as follows: 
 \begin{eqnarray*}
 {\bf H_1^{\it n}}&:& \mbox{ With probability 1, } \lambda=\lambda_0+d_n \Delta, \mbox{ where }  \lambda_0\, :\, [0,T]\to \mathbb R_+ \mbox{ is a bounded}\\
 & & \mbox{ non-null function and } \Delta=(\Delta_t)_{t\in [0,T]} \mbox{ is a stochastic process with }\\
 & &  \sup_{t\in [0,T]} \mathbb E \Delta_t^6 <\infty \mbox{ and}\ {\rm var} \big(\int_0^t \Delta_s {\rm d}s\big)>0 \mbox{ for some } t\in [0,T].
  \end{eqnarray*} 
 When $n$ gets larger and ${\bf H_1^{\it n}}$ holds, $N$ is becomes closer to a Poisson process. Thus, $(d_n)_n$ has to be understood as a separation speed from ${\bf H_1^{\it n}}$ to the null hypothesis ${\bf H_0}$. 
 
 In particular, next result states that in view of a consistent test for ${\bf H_0}$ vs ${\bf H_1^{\it n}}$, it is necessary and sufficient that $d_n^2$ tends to 0 slower than $1/\sqrt{n}$. 
   
\begin{theo} \label{local} Let $B=(B_t)_{t\in\mathbb R_+}$ be a standard Brownian Motion on the real line. Assume that ${\bf H_1^{\it n}}$ holds, and denote by $m_0$ and $v$ the functions defined for all $t\in [0,T]$ by
$$m_0(t)=\int_0^t \lambda_0(s){\rm d}s, \mbox{ and } v(t)={\rm var} \big(\int_0^t \Delta_s {\rm d}s\big).$$
Moreover, we let $I_0^2=\int_0^T (T-t)m_0(t)^2 {\rm d}t.$

\smallskip
(i) If $\sqrt{n} \, d_n^2\to \infty$, then 
$$\sqrt{n}\frac{\hat S_1}{\hat m(T)}\, \stackrel{\rm prob.}\longrightarrow \, +\infty, \mbox{ and } \sqrt{n}\frac{\hat S_2}{\hat I}\, \stackrel{\rm prob.}\longrightarrow \, +\infty.$$

(ii) If $\sqrt{n} \, d_n^2\to d<\infty$, then
\begin{eqnarray*}
& & \sqrt{n} \frac{\hat S_1}{\hat m(T)} \, \stackrel{\rm{(law)}}\longrightarrow  \, \frac{1}{m_0(T)} \sup_{t\le T} \big(B_{2m_0(t)^2}+dv(t)\big), \mbox { and}\\
&  & \sqrt{n} \frac{\hat S_2}{\hat I} \, \stackrel{\rm{(law)}}\longrightarrow  \, 2\mathcal N(0,1)+\frac{d}{I_0} \int_0^T v(t){\rm d}t.
\end{eqnarray*}
\end{theo} 

In the problem ${\bf H_0}$ vs ${\bf H_1^{\it n}}$, we consider the tests defined by the rejection regions in \eqref{RR}, with $\alpha\in ]0,1[$. For a power study, we assume from now on that ${\bf H_1^{\it n}}$ holds. By above, if $\sqrt{n} d_n^2\to d<\infty$, 
\begin{eqnarray}
\label{r2} 
 \lim_{n\to \infty} \mathbb P\big(R_2(\alpha)\big) = 1-\Phi\Big(q_{1-\alpha}-\frac{d}{2I_0} \int_0^T v(t){\rm d}t\Big)<1,
 \end{eqnarray}
where $\Phi$ stands for the repartition function of the $\mathcal N(0,1)$ distribution. However, we only have 
\begin{eqnarray*}
\limsup_{n\to \infty} \mathbb P \big(R_1(\alpha)\big)\le \mathbb P \Big(\frac{1}{m_0(T)} \sup_{t\le T} \big(B_{2m_0(t)^2}+dv(t)\big)\ge \sqrt{2} q_{1-\alpha/2}\Big)<1,
\end{eqnarray*}
according to the Portmanteau Theorem, as the limit distribution may have a mass at point $\sqrt{2} q_{1-\alpha/2}$. At least, we deduce from above and part $(i)$ of Theorem \ref{local} that both tests defined by $R_1(\alpha)$ and $R_2(\alpha)$ are consistent if, and only if $\sqrt{n}d_n^2\to\infty$. 

In the rest of the section, we assume that $\sqrt{n}d_n^2\to d<\infty$. For a comparison of the tests, we need an additional assumption ensuring that the limit distribution associated with statistic $\hat S_1$ is continuous. To this aim, we suppose now that $\lambda_0(t)>0$ for all $t\in [0,T]$, and we let $\ell_0$ the function such that $\ell_0(t)=2m_0(t)^2$ for all $t\in [0,T]$.  Then, $\ell_0$ is a continuous and increasing function with $\ell_0(0)=0$, and 
\begin{eqnarray}
\label{sup0} 
\sup_{t\le T} \big(B_{\ell_0(t)}+dv(t)\big)=\sup_{s\le \ell_0(T)}\big(B_s+dv\circ \ell_0^{-1}(s)\big).
\end{eqnarray}
Now observe that function $v\circ \ell_0^{-1}$ is absolutely continuous. Thus, by the Girsanov Theorem (Revuz and Yor, 1999), there exists a probability measure $Q$ such that the distribution of the random variable in \eqref{sup0}  equals the distribution under $Q$ of the supremum of a standard Brownian Motion over $[0,\ell_0(T)]$. According to Proposition III.3.7 in the book by Revuz and Yor (1999), this distribution is continuous, which proves that the distribution of the random variable in \eqref{sup0} is also continuous. As consequence, the Portmanteau Theorem and Theorem \ref{local} $(ii)$ give  
$$\lim_{n\to \infty} \mathbb P \big(R_1(\alpha)\big)= \mathbb P \Big(\frac{1}{m_0(T)} \sup_{s\le \ell_0(T)} \big(B_s+dv\circ \ell_0^{-1}(s)\big)\ge \sqrt{2} q_{1-\alpha/2}\Big).$$
Unfortunately, the latter probability is not known, except for some special cases that we now study. 

\begin{figure}[H]
\begin{center}
\includegraphics[width=6.8cm,height=7cm]{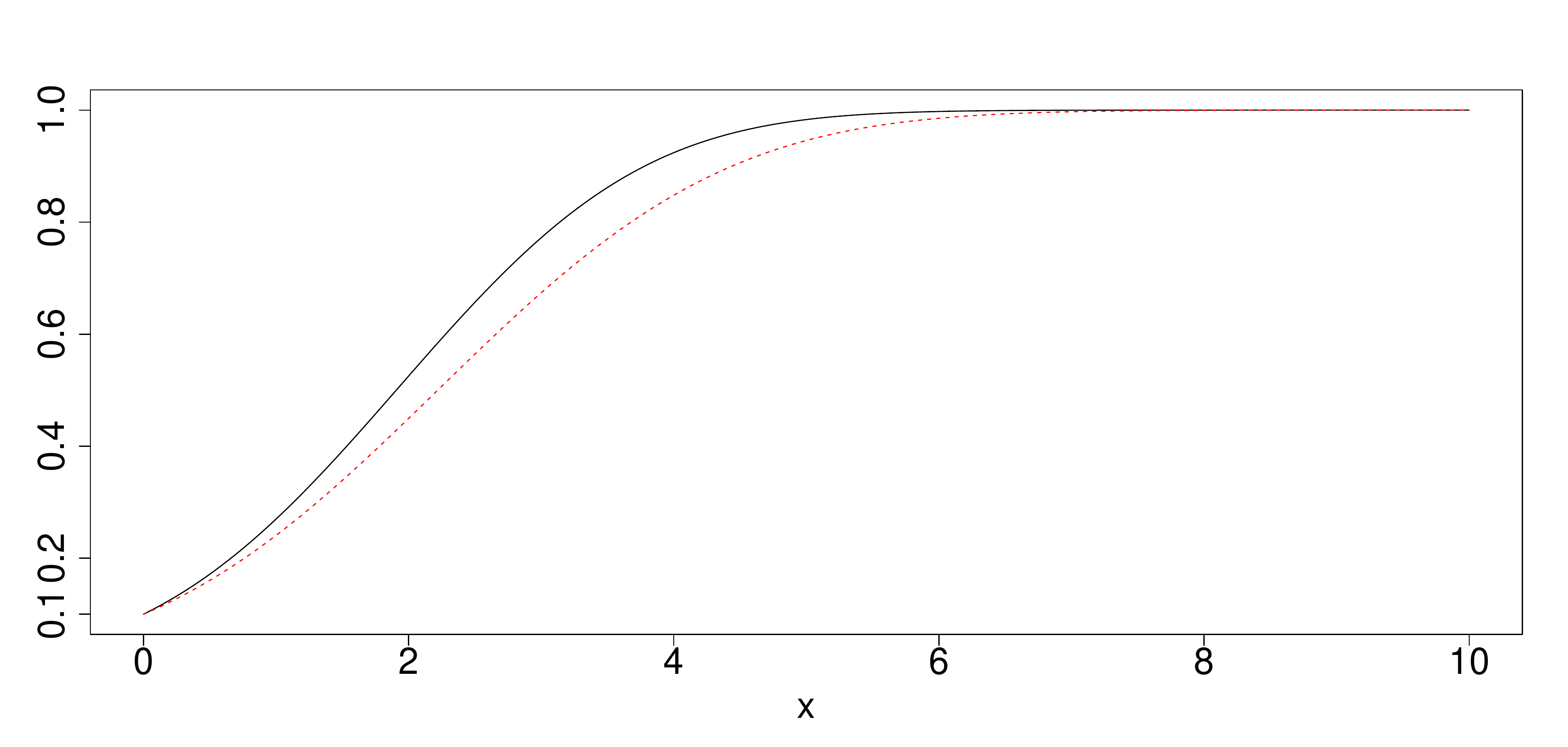}
\includegraphics[width=6.8cm,height=7cm]{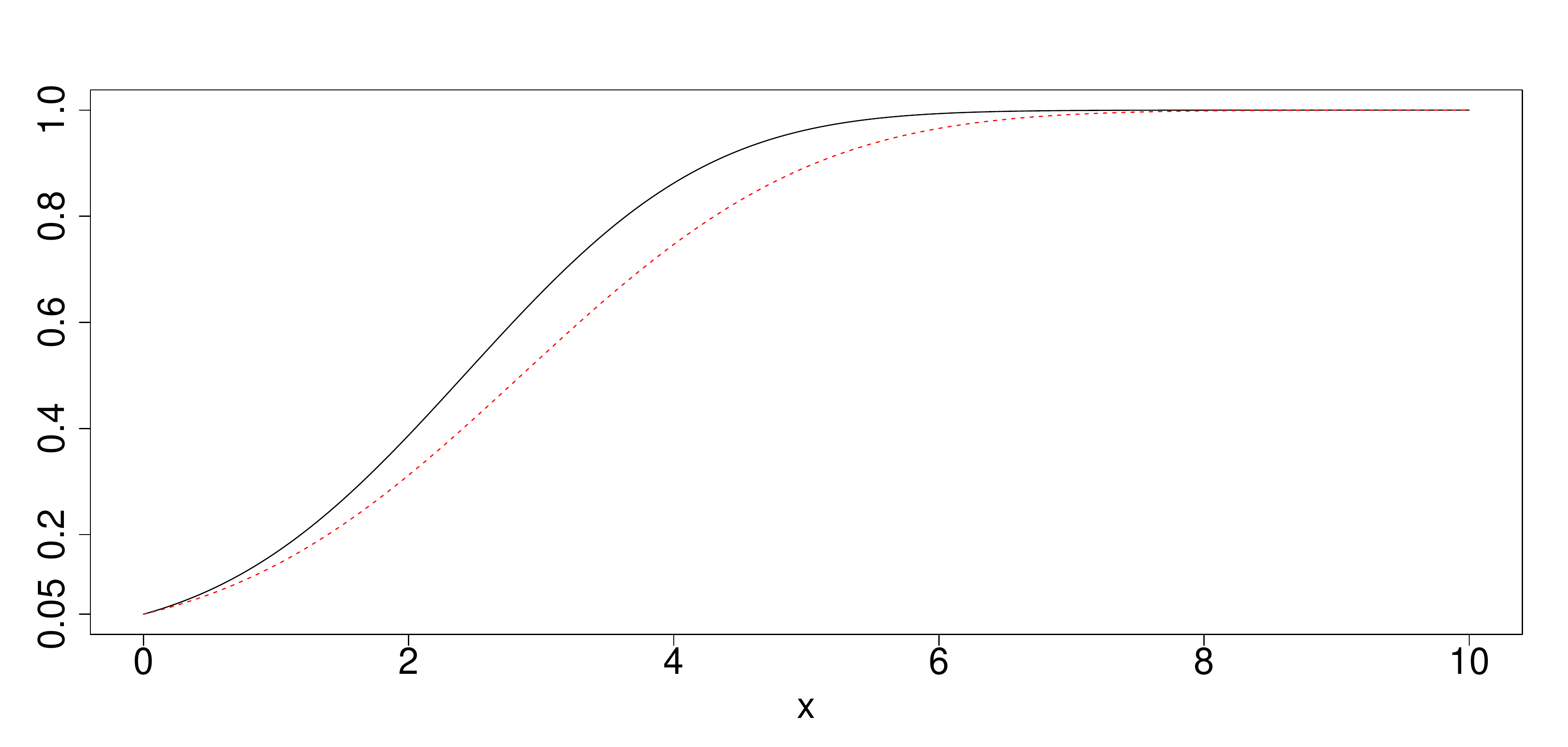}
\end{center}
\caption{Plots of $g_1$ (full lines) and $g_2$ (dashed lines), with $\alpha=0.1$ for the curves on the left, $\alpha=0.05$ for the curves on the right.\label{comp}}
\end{figure}

In addition to ${\bf H_1^{\it n}}$, we suppose that $\lambda_0(t)=\lambda_0>0$ and $\Delta_t=Z$ for all $t\in [0,T]$, where $Z$ is a random variable with variance $w^2$. In particular, we study the case of a small deviation from an homogeneous Poisson process. Then, formula \eqref{sup0} writes 
$$\sup_{t\le T} \big(B_{\ell_0(t)}+dv(t)\big)=\sup_{s\le 2\lambda_0^2 T^2}\big(B_s+\frac{dw^2}{2\lambda_0^2}s\big).$$
Obviously, we have in this case $v(t)=w^2t^2$, $m_0(t)=\lambda_0 t$ and $I_0^2=\lambda_0^2 T^4/12$ for all $t\in [0,T]$. Setting $x=dw^2T$, we obtain with Theorem \ref{local} $(ii)$ and the distribution of the supremum of a drifted Brownian Motion (see p. 250 in the book by Borodin and Salminen, 2002):
\begin{eqnarray*}
\lim_{n\rightarrow \infty}\mathbb{P}\big(R_1(\alpha)\big)&=&\exp\big(\frac{\sqrt{2}}{\lambda_0}x q_{1-\alpha/2}\big)\Big(1-\Phi\big(q_{1-\alpha/2}+\frac{x}{\sqrt{2}\lambda_0}\big)\Big)\\
& & +1-\Phi\big(q_{1-\alpha/2}-\frac{x}{\sqrt{2}\lambda_0}\big).
\end{eqnarray*}
Moreover, by \eqref{r2}: 
$$\lim_{n\rightarrow \infty}\mathbb{P}\big(R_2(\alpha)\big)=1-\Phi\big(q_{1-\alpha}-\frac{x}{\sqrt{3}\lambda_0}\big).$$
Assume $\lambda_0=1$ and denote by $g_1$ and $g_2$ respectively the above functions of $x$. Figure \ref{comp} shows a comparison of these two quantities as functions of $x$ and for $\alpha=0.05$ or $\alpha=0.1$.
This example suggests a better power for the test induced by $R_1(\alpha)$. The numerical experiments given in the next section also suggest a better power for the first test.

\section{Simulation study}\label{section:simu}

In this section we illustrate the good properties of our tests with a simulation study. We consider $n_{MC}$ replications of Monte Carlo simulations and we study the performances of our tests in terms of asymptotic level and rejection power.

Throughout this section, we fix $T=1$. In the following, function $\lambda$ denotes the intensity function of $N$, \textit{i.e.} the first derivative of its cumulative intensity function $\Lambda$.

\subsection{Level study}
We consider the following model for the asymptotic level study,
\begin{align}
\lambda(t)=\beta t^{\beta-1},\text{ with }\beta>0. \label{eq:level}
\end{align}
This is the intensity function of a so-called Weibull process, which is frequently used in Reliability Theory for instance. Remark that function $\lambda$ is decreasing for $\beta<1$, constant for $\beta=1$ and increasing for $\beta>1$.

In Table~\ref{tab:level} we evaluate the empirical rejection frequency of both tests using the rejection regions defined in \eqref{RR} for levels $\alpha=5\%$ and $\alpha=10\%$. This evaluation is based on $10,000$ Monte Carlo simulations (for each value of $\beta\in\{1/2,1,2\}$), with $n=100$ and $500$ replications of the model. We note that test statistic $\hat S_1$ has a similar behavior than test statistic $\hat S_2$ at this range. For both statistics, the empirical rejection frequency is close to the nominal value even with the smallest sample size $n=100$. 
%

\begin{table}[ht]
\begin{center}
\begin{tabular}{c|cccccccccc}
\multicolumn{2}{l}{}&\multicolumn{2}{c}{$\beta=1/2$}&&\multicolumn{2}{c}{$\beta=1$}&&\multicolumn{2}{c}{$\beta=2$}&\\
\cline{3-4} \cline{6-7} \cline{9-10}
\multicolumn{2}{l}{$n$}&$100$&$500$&&$100$&$500$&&$100$&$500$&\\
\hline
&&\multicolumn{9}{c}{$\alpha=5\%$}\\
\cline{2-11}\\
$\hat S_1$&&$6.55$&$5.84$&&$5.99$&$5.52$&&$5.69$&$5.56$&\\
$\hat S_2$&&$6.27$&$5.38$&&$5.74$&$5.44$&&$6.00$&$5.71$&\\
\\
\hline
&&\multicolumn{9}{c}{$\alpha=10\%$}\\
\cline{2-11}\\
$\hat S_1$ && $10.96$ & $10.29$ && $10.25$ & $9.95$ && $10.45$ & $10.25$ & \\
$\hat S_2$ && $10.55$ & $10.05$ && $10.38$ & $10.02$ && $10.29$ & $10.06$ & \\
\end{tabular}
\caption{Empirical rejection frequency ($\%$) for $n_{MC}=10,000$ \label{tab:level}.}
\end{center}
\end{table}

\subsection{Rejection power study}

\begin{figure}[ht]
\begin{center}
\includegraphics[scale=0.9]{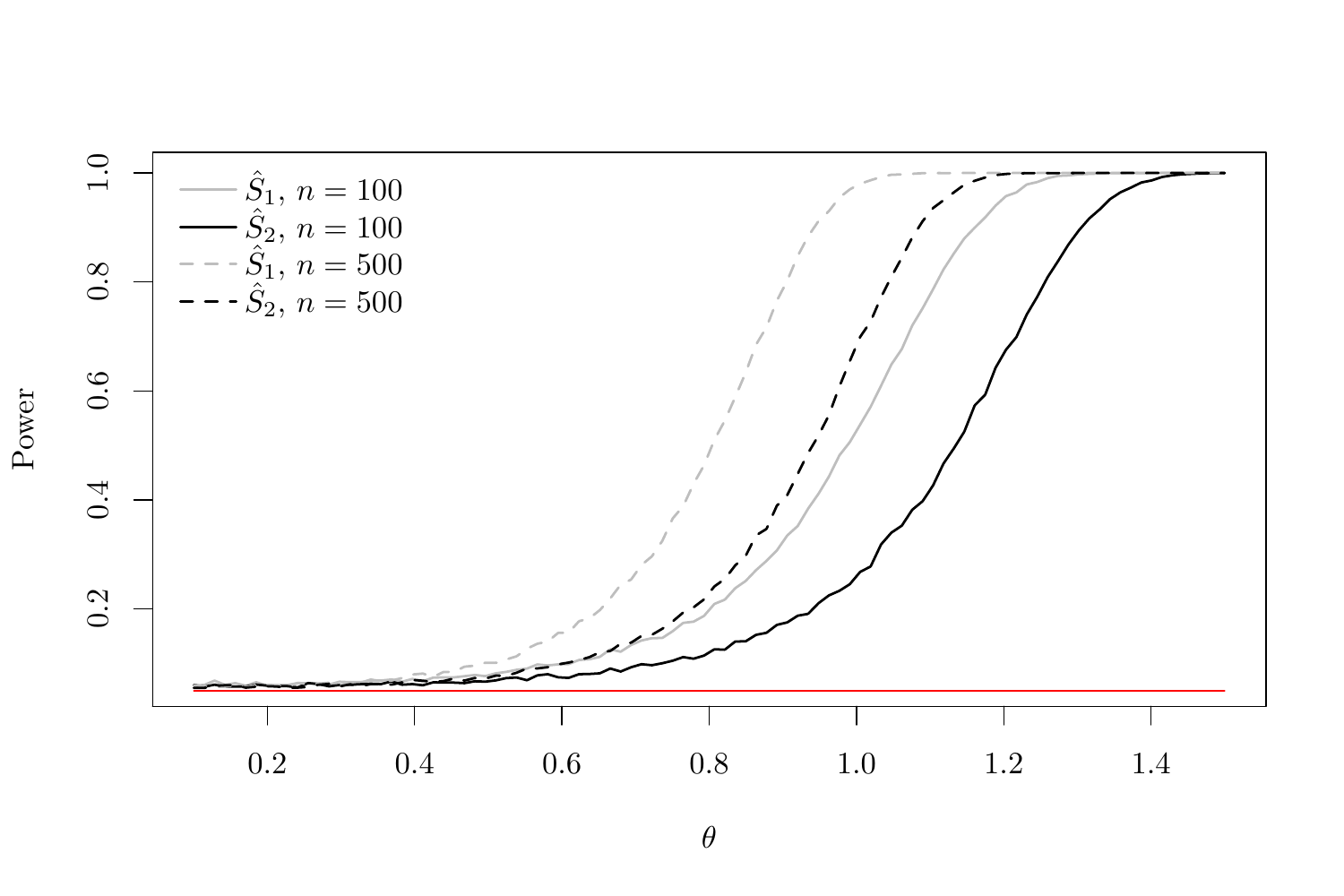}
\end{center}
\caption{\label{fig:power} Empirical rejection frequency under $\eqref{eq:power}$ for $n_{MC}=10,000$, $100$ and $500$ trajectories and $\alpha=5\%$. The red horizontal line represents the value of $\alpha$.}
\end{figure}

\noindent {\it Model 1.}Ê For the asymptotic power study, we first consider the model defined by 
\begin{align}
\lambda(t)=\exp(\theta Zt), \label{eq:power}
\end{align}
with $Z\sim2+\mathrm{Beta}(\frac{1}{2},\frac{1}{2})$ and $\theta\in [0,1]$. 
Figure~\ref{fig:power} represents the empirical rejection frequency for different values of $\theta$ between $0$ and $1$ in model \eqref{eq:power}. For $\theta=0$, the simulated process is an homogeneous Poisson process with intensity $1$. The process deviates from the homogeneous Poisson process with the increase of $\theta$. We observe on figure~\ref{fig:power} that both tests catch this behavior for small $\theta$'s. In both cases the power goes to $1$ for higher values of $\theta$. The test statistic $\hat S_1$ has better performances as its power curve increases faster that of test statistic $\hat S_2$.\\

\begin{figure}[ht]
\begin{center}
\includegraphics[scale=0.9]{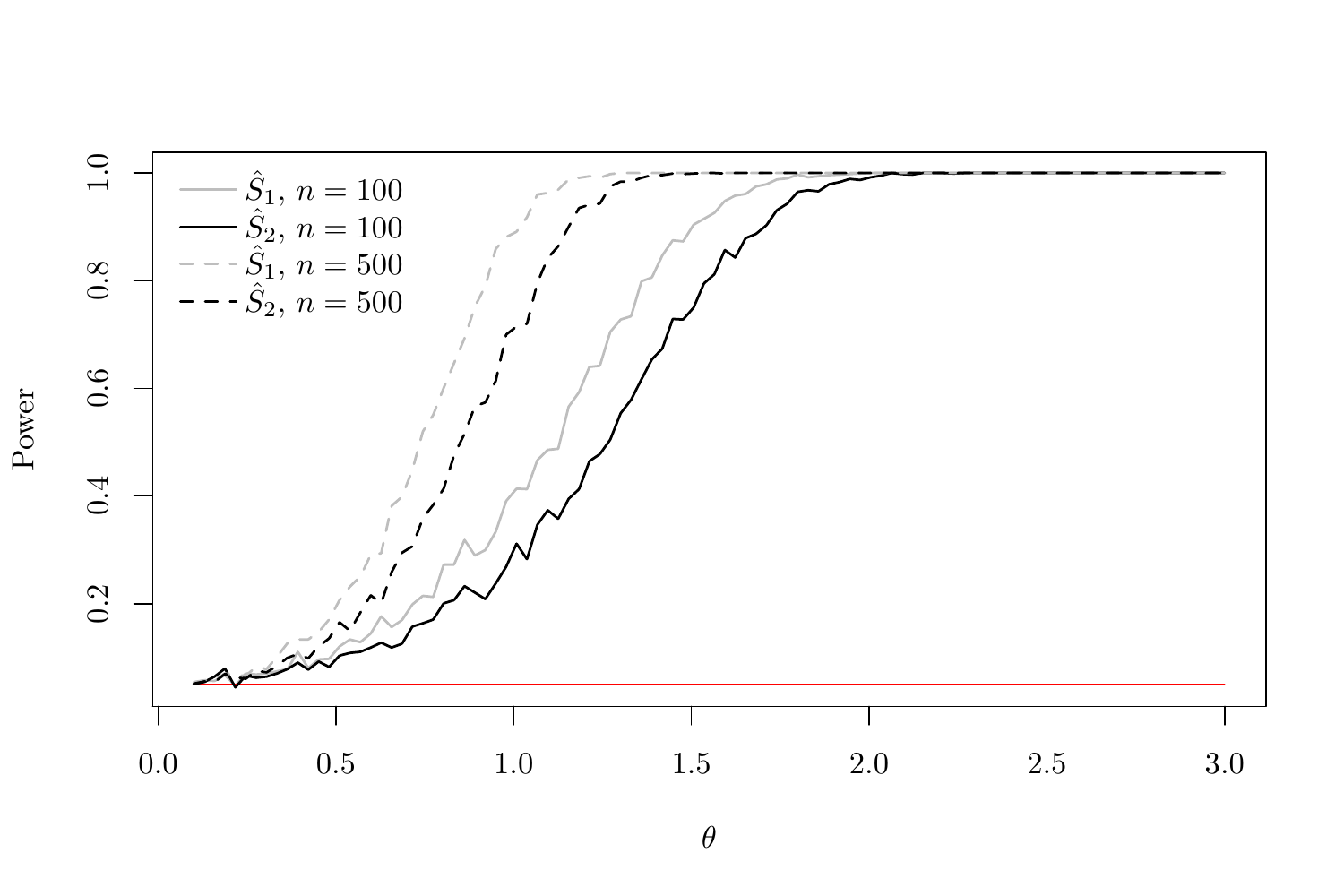}
\end{center}
\caption{\label{fig:powerZt} Empirical rejection frequency under $\eqref{eq:powerZt}$ for $n_{MC}=1,000$, $100$ and $500$ trajectories and $\alpha=5\%$. The red horizontal line represents the value of $\alpha$.}
\end{figure}

\noindent {\it Model 2.} We define the second model as follows,
\begin{align}
\lambda(t)=\exp(\theta \sin(Z_t)), \label{eq:powerZt}
\end{align}
with $(Z_t)_{t\in[0,1]}$ a standard Brownian Motion and $\theta\in [0,1]$. This model differs from the previous one as the covariate depends on the time variable. Figure~\ref{fig:powerZt} represents the empirical rejection frequency for $\theta$ varying in $[0,3]$ in model \eqref{eq:powerZt}. The same remarks as for model~\eqref{eq:power} apply here. Note that the two test statistics look closer on this example.

\section{Application to real data}\label{section:data}
\subsection{Analysis of some arrival times in a call center}
The use of Poisson processes has been often considered as a first approach for modeling the arrival times in call centers and more generally queue systems, e.g. see Asmussen (2003) for the nice theory developed around this assumption.

As in the papers by Brown et al. (2005) and Mandelbaum et al. (2000), we consider a call center for an anonymous Israel's banks. A description of the calls received over the year $1999$ is available online\footnote{http://ie.technion.ac.il/serveng/callcenterdata/}. The call center is open on weekdays (Sunday to Thursday in Israel), from 7\AM to midnight, and the calls are separated in different classes, depending on the needs of the customers. Each call can be described as follows. A customer calls one of the phone numbers of the call center. Except for rare busy signals, the customer is then connected to an interactive voice respond unit (IVR or VRU) and identifies himself/herself. While using the VRU, the customer receives recorded information. He/She can choose to perform some self-service transactions ($\sim65\%$ of the calls) or indicate the need to speak with an agent ($\sim35\%$). Here, we are interested in the latter, which represents roughly $30,000$ to $40,000$ calls per month.
Each call record in the database also includes a categorical description of the type of service requested.
The main call types are regular (PS in the database),
stock transaction (NE), new or potential customer (NW), and Internet assistance (IN). Mandelbaum et al. (2000) and Brown et al. (2005) described the process of collecting and cleaning the data and provided complete descriptive analysis of the data.

In this study, we concentrate on IN calls arriving between 3:25\PM and 3:35\PM on all weekdays of year 1999. Times at which calls enter the VRU represent the arrival times of a counting process. The dataset then consists in 258 trajectories of this type, that we can assume to be realizations of {\it i.i.d.} Cox processes.  

The results of the statistical study are presented in Table~\ref{table:callcenter}.
One can see that the null hypothesis ${\bf H_0}$ is highly rejected using both statistics.
This result suggests that even on a short period of time, these arrival times, which depend on a complex human behavior, seem to be strongly influenced by some covariates. One might easily imagine that weather conditions or other company intrinsic variables ({\it e.g.} number of recent opened accounts) could reduce this overdispersion and help to explain the number of IN phone calls. 

\begin{table}[ht]
\centering
\begin{tabular}{c|cc}
&$\hat S_1$&$\hat S_2$\\
\hline
$p$-values&$1.95\times10^{-6}$&$1.05\times10^{-6}$
\end{tabular}
\caption{$p$-values of both tests for the call center dataset.\label{table:callcenter}}
\end{table}

Finally mention that Brown et al. (2005) also studied IN calls but did not reject the Poisson assumption. However, their study consists in testing the exponential distribution for the interarrival times of IN calls occurring in a single day, which is not compatible with our asymptotic and cannot help to determine if some covariates influence the daily calls process.

\subsection{Analysis of the scoring times of a soccer team}
As seen in Heuer et al. (2010), Poisson processes may also be used to model scoring goals during a soccer match. Nevertheless, one could suspect the influence of some covariates such as the behavior of the spectators or fitness fluctuations of the team under study. 
Thus, we propose to test the Poisson process assumption ${\bf H_0}$ for the scoring times of Arsenal soccer club first team.

To this end, we collected on the SoccerSTATS.com website their scoring times (in minutes) for each match in "Premier League" over three seasons (from 2012 to 2015), for a total of $229$ matches. 
For each match, the scoring times of the team define the jump times of 229 counting paths. 
We can assume that these data are {\it i.i.d.} realizations of Cox processes. 

The results of the statistical study are presented in Table~\ref{table:soccer}. 
For both statistics $\hat S_1$ and $\hat S_2$, we cannot reject the null hypothesis ${\bf H_0}$ and the Poisson process seems to be a reasonable 
approximation for these counting processes.  
\begin{table}[ht]
\centering
\begin{tabular}{c|cc}
&$\hat S_1$&$\hat S_2$\\
\hline
$p$-values&$0.419$&$0.298$
\end{tabular}
\caption{$p$-values of both tests for the soccer goals dataset.\label{table:soccer}}
\end{table}

Recall that the analysis in Heuer et al. (2010) also suggests that the Poisson process is relevant for modeling the scoring goals of a given team in the German Bundesliga.

\section{Proof of Theorem \ref{lim1}}
\label{section:proofs}

In the rest of the paper, we assume for notational simplicity that $T=1$. We let $\hat M=(\hat M_t)_{t\in [0,1]}$ the process defined by $\hat M_t=\hat \sigma^2(t)-\hat m(t)$ and $\hat \tau=(\hat \tau_t)_{t\in [0,1]}$ is the process such that for all $t\in [0,1]$, $$\hat \tau_t=\frac{4}{n-1} \int_0^t \hat \sigma^2(s){\rm d}m(s).$$
Martingale properties and predictable $\sigma$-field are implicitly with respect to the natural filtration generated by the sample $N^{(1)},\cdots,N^{(n)}$. As usual, $\langle X \rangle$ stands for the predictable quadratic variation of the martingale $X$.

\subsection{Auxiliary results} 

\begin{lem} \label{compensateur} Under ${\bf H_0}$, the process $\hat M$ is a martingale, and $\langle \hat M\rangle=\hat \tau$. \end{lem} 

\noindent {\bf Proof.} First we prove that $\hat M$ is a martingale. Observe that for all $t\in [0,1]$, 
\begin{equation}
\label{eq:s}
\hat \sigma^2(t)=\frac{1}{n-1} \sum_{k=1}^n \big(N_t^{(k)}\big)^2-\frac{n}{n-1} \hat m(t)^2.
\end{equation}
We now compute $\hat M_t$ as a sum of stochastic integrals. In the sequel, we let $\bar N^{(k)}=N^{(k)}-m$. Note that $\bar N^{(k)}$ is a martingale. According to the integration by parts formula (Proposition 0.4.5 in the book by Revuz and Yor, 1999), we have 
\begin{eqnarray}
\sum_{k=1}^n \big(N_t^{(k)}\big)^2 & = & \sum_{k=1}^n \Big[2\int_0^t N_{s^-}^{(k)}{\rm d}N_s^{(k)}+N^{(k)}_t\Big]\nonumber\\
& = & 2\sum_{k=1}^n \int_0^t N_{s^-}^{(k)}{\rm d}N_s^{(k)}+n\hat m(t)\nonumber \\
& = & 2\sum_{k=1}^n \int_0^t N_{s^-}^{(k)}{\rm d}\bar N^{(k)}_s+2n\int_0^t \hat m(s^-){\rm d}m(s)+n\hat m(t). \label{eq:1} 
\end{eqnarray}
Moreover, by the integration by parts formula, 
$$\hat m(t)^2=2\int_0^t \hat m(s^-){\rm d}\hat m(s)+\sum_{s\le t} \big(\Delta \hat m(s)\big)^2.$$
Using the fact that two independent Poisson processes do not jump at the same time (Proposition XII.1.5 in the book by Revuz and Yor, 1999), we deduce that 
\begin{eqnarray*}
\sum_{s\le t} \big(\Delta \hat m(s)\big)^2 & = & \frac{1}{n^2} \sum_{s\le t} \Big(\sum_{k=1}^n \Delta N_s^{(k)}\Big)^2= \frac{1}{n^2} \sum_{s\le t} \sum_{k=1}^n \Delta N_s^{(k)}\\
& = & \frac{1}{n} \hat m(t).
\end{eqnarray*}
Hence, 
$$\hat m(t)^2=2\int_0^t \hat m(s^-){\rm d}\hat m(s)+\frac{1}{n}\hat m(t).$$
Then, combining \eqref{eq:s} and \eqref{eq:1}, we obtain 
\begin{eqnarray}
\label{final}
\hat M_t &= & -\hat m(t)+\frac{1}{n-1} \sum_{k=1}^n \big(N_t^{(k)}\big)^2-\frac{n}{n-1}\hat m(t)^2\nonumber \\
& = & \frac{2}{n-1}\sum_{k=1}^n \int_0^t N_{s^-}^{(k)}{\rm d}\bar N_s^{(k)}-\frac{2n}{n-1}\int_0^t \hat m(s^-){\rm d}\big(\hat m(s)-m(s)\big)\label{fin1}
\end{eqnarray}
Since $\hat m-m$ and each of the $\bar N^{(k)}$'s are martingales and the integrands are predictables, we deduce that all integrals in this formula are local martingales. It is a classical exercise to prove that  they are of class DL (see Definition IV.1.6 in the book by Revuz and Yor, 1999), so that they are martingales, as well as $\hat M$. 

In view of computing the predictable quadratic variation of $\hat M$, we first observe that by the integration by parts formula, 
\begin{equation} 
\label{eq:11}
\hat M^2_t=2\int_0^t \hat M_{s^-}{\rm d}\hat M_s+\sum_{s\le t} \big(\Delta \hat M_s\big)^2.
\end{equation} 
But, by \eqref{fin1}, 
$$\Delta \hat M_s=\frac{2}{n-1} \sum_{k=1}^n N_{s^-}^{(k)}\Delta N_s^{(k)}-\frac{2n}{n-1} \hat m(s^-)\Delta \hat m(s).$$
Again, we shall make use of the fact that two Poisson processes do not jump at the same time. Hence, if $s$ is a time-jump for $N^{(k)}$, 
\begin{equation}\label{saut} \Delta \hat M_s = \frac{2}{n-1} \big(N^{(k)}_{s^-}-\hat m(s^-)\big)= \frac{2}{n-1} \big(N_{s^-}^{(k)}-\hat m(s^-)\big) \Delta N_s^{(k)}.\end{equation}
Thus, 
\begin{eqnarray*}
\sum_{s\le t} \big(\Delta \hat M_s\big)^2 & = & \sum_{s\le t} \sum_{k=1}^n \big(\Delta \hat M_s\big)^2\mathbf 1_{\{\Delta N_s^{(k)}=1\}}\\
& = &  \frac{4}{(n-1)^2}\sum_{k=1}^n \int_0^t \big(N_{s^-}^{(k)}-\hat m(s^-)\big)^2{\rm d}N_s^{(k)}.
\end{eqnarray*}
By \eqref{eq:11}, we have 
\begin{eqnarray*}
\hat M_t^2 & = & 2\int_0^t \hat M_{s^-}{\rm d}\hat M_s+\frac{4}{(n-1)^2}\sum_{k=1}^n \int_0^t \big(N_{s^-}^{(k)}-\hat m(s^-)\big)^2{\rm d}\bar N_s^{(k)}\\
& & + \frac{4}{(n-1)^2}\sum_{k=1}^n \int_0^t \big(N_{s^-}^{(k)}-\hat m(s^-)\big)^2{\rm d}m(s).
\end{eqnarray*}
As above, we can conclude from the fact that both $\hat M$ and $\bar N^{(k)}$ are martingales that the first two terms on the right-hand side are martingales. Last term, namely 
$$\frac{4}{n-1} \int_0^t \hat \sigma^2(s^-){\rm d}m(s)=\frac{4}{n-1} \int_0^t \hat \sigma^2(s){\rm d}m(s),$$
where equality holds by continuity of $m$, is predictable. Hence, it is the predictable quadratic variation of $\hat M$. $\Box$

\begin{lem} \label{lemma:saut} Under ${\bf H_0}$, we have :
$$n\mathbb E \sup_{t\le 1} |\Delta \hat M_t|^2\to 0.$$
\end{lem} 

\noindent {\bf Proof.} First observe that for all $u\in ]0,1[$ : 
\begin{eqnarray}
n\mathbb E \sup_{t\le 1} |\Delta \hat M_t|^2 & \le & u+n\mathbb E \sup_{t\le 1} |\Delta \hat M_t|^2\mathbf 1_{\{\sqrt{n} \sup_{t\le 1} |\Delta \hat M_t|> u\}}\nonumber\\
& \le & u+2\int_0^\infty x \mathbb P\big(\sqrt{n} \sup_{t\le 1} |\Delta \hat M_t|\mathbf 1_{\{\sqrt{n} \sup_{t\le 1} |\Delta \hat M_t|>u\}}>x\big) {\rm d}x\nonumber \\
& \le & 2u + 2\int_u^\infty x \mathbb P(\sqrt{n} \sup_{t\le 1} |\Delta \hat M_t| >x){\rm d}x. \label{doob}
\end{eqnarray}
But, according to \eqref{saut}, if $t$ is a time-jump for $N^{(k)}$, we have 
$$\Delta \hat M_t=\frac{2}{n-1}\big(N_{t^-}^{(k)}-\hat m(t^-)\big),$$
and hence, 
$$\sup_{t\le 1}|\Delta \hat M_t|\le \frac{2}{n-1} \sup_{k\le n} \sup_{t\le 1} |N_{t}^{(k)}-\hat m(t)|.$$
Thus, for all $x>0$ : 
\begin{eqnarray*}
\mathbb P(\sqrt{n}\sup_{t\le 1} |\Delta \hat M_t|\ge x)& \le & n \mathbb P\Big(\sup_{t\le 1} |N_t^{(1)}-\hat m(t)|\ge \frac{(n-1)x}{2\sqrt{n}}\Big)\\
& \le & 8 \sup_{t\le 1} \mathbb E |N_t^{(1)}-\hat m(t)|^3 \frac{n^{5/2}}{(n-1)^3 x^2},
\end{eqnarray*}
where the last inequality is due to Doob's Inequality (see Revuz and Yor, 1999) applied to the martingale $N^{(1)}-\hat m$. A direct calculation shows that that exists a constant $C>0$ (independent of $n$) such that 
$$\sup_{t\le 1} \mathbb E |N_t^{(1)}-\hat m(t)|^3 \le C.$$
Consequently, by \eqref{doob} : 
\begin{eqnarray*} 
n\mathbb E \sup_{t\le 1} |\Delta \hat M_t|^2& \le &  2u+ \frac{16Cn^{5/2}}{(n-1)^3} \int_u^\infty \frac{{\rm d}x}{x^2}\\
& \le & 2u+ \frac{16C n^{5/2}}{u(n-1)^3}.
\end{eqnarray*}
Taking for instance $u=n^{-1/4}$ gives the result. $\Box$

\subsection{Proof of Theorem \ref{lim1}} According to Theorem VIII.3.22 in the book by Jacod and Shiryaev (2003), the sequence of square integrable martingales $(\sqrt{n}\hat M)_n$ converges in distribution to a continuous Gaussian martingale $M$ such that $\langle M \rangle=2m^2$ if, for all $t\in [0,1]$ and $\varepsilon >0$, 
\begin{equation}
\label{proprietes} \langle \sqrt{n} \hat M\rangle_t\to 2m(t)^2 \mbox{ and } \int_{\mathbb R \times [0,t]} |x|^2 \mathbf 1_{\{|x|>\varepsilon\}} \nu_n({\rm d}x,{\rm d}s)\to 0,
\end{equation}
both in probability, where $\nu_n$ stands for the predictable compensator of the random jump measure associated to the martingale $\sqrt{n} \hat M$. Regarding the first property, we know from Lemma \ref{compensateur} that 
$$\lim_{n\to\infty} \langle \sqrt{n} \hat M\rangle_t=\lim_{n\to \infty} n\hat \tau_t=4\int_0^t \sigma^2(s){\rm d}m(s),$$
in probability. Since $\sigma^2=m$ under ${\bf H_0}$, we deduce that 
$\lim_{n\to\infty} \langle \sqrt{n} \hat M\rangle_t=2m(t)^2$
in probability. In order to prove the second property in \eqref{proprietes}, we fix $\varepsilon >0$ and we let $U$ and $V$ be the processes defined for all $t\in [0,1]$ by 
$$U_t =  \int_{\mathbb R \times [0,t]} |x|^2 \mathbf 1_{\{|x|>\varepsilon\}} \nu_n({\rm d}x,{\rm d}s) \mbox{ and }
V_t = n\sum_{s\le t} |\Delta \hat M_s|^2 \mathbf 1_{\{\sqrt{n} |\Delta \hat M_s|>\varepsilon\}}.
$$
Observing that $U$ is L-dominated by the increasing adapted process $V$, we deduce from the Lenglart Inequality (see p.35 in the book by Jacod and Shiryaev, 2003) that for all $t\in [0,1]$ and $\alpha,\eta >0$ : 
$$\mathbb P(U_t\ge \eta) \le \frac{1}{\eta} \big(\alpha+\mathbb E \sup_{s\le t} \Delta V_s\big)+\mathbb P(V_t\ge \alpha).$$
But, $\{V_t>0\}=\{\sqrt{n} \sup_{s\le t} |\Delta \hat M_s|>\varepsilon\}$ and $\sup_{s\le t} \Delta V_s\le n\sup_{s\le t} |\Delta \hat M_s|^2$. Thus, letting $\alpha\searrow 0$, we obtain with the help of Markov's Inequality :
\begin{eqnarray*}
\mathbb P(U_t\ge \eta) &\le &  \frac{n}{\eta} \mathbb E \sup_{s\le t} |\Delta \hat M_s|^2+\mathbb P(\sqrt{n} \sup_{s\le t} |\Delta \hat M_s|>\varepsilon)\\
& \le & \Big(\frac{1}{\eta}+\frac{1}{\varepsilon^2}\Big) n \mathbb E \sup_{s\le t} |\Delta \hat M_s|^2.
\end{eqnarray*}
We conclude from Lemma \ref{lemma:saut} that $U_t$ converges to 0 in probability. Hence, both properties in \eqref{proprietes} are satisfied so that the sequence of square integrable martingales $(\sqrt{n}\hat M)_n$ converges in distribution to a continuous Gaussian martingale $M$, such that $\langle M \rangle=2m^2$. The Dambis-Dubins-Schwarz Theorem (see Theorem V.1.6 in the book by Revuz and Yor, 1999) then gives $M=B_{2m^2}$, where $B$ is a standard real Brownian Motion. $\Box$ 

\section{Proof of corollary \ref{limite}} 

 {\it (i)} Let $\mathbb D$ be the space of c\`adl\`ag functions from $[0,1]$ to $\mathbb R$, equipped with the Skorokhod topology. By continuity of the application $\mathbb D \ni x\mapsto \sup_{t\le T} x(t)$, we deduce from Theorem \ref{lim1} that
$$\sqrt{n}\, \hat S_1=\sqrt{n} \sup_{t\le 1}\big(\hat \sigma^2(t)-\hat m(t)\big)\stackrel{{\rm (law)}}\longrightarrow \sup_{t\le 1} B_{2m(t)^2}=\sup_{t\le 2m(1)^2} B_t.$$
According to the reflection principle (Proposition III.3.7 in the book by Revuz and Yor, 1999), the distribution of the latter term is $\sqrt{2}\, m(1) |\mathcal N(0,1)|$, hence the result with $\hat S_1$. Similarly, by continuity of $\mathbb D\ni x\mapsto \int_0^1 x(t){\rm d}t$, we have 
$$\sqrt{n} \hat S_2=\sqrt{n} \int_0^1 \big(\hat \sigma^2(t)-\hat m(t)\big){\rm d}t \stackrel{{\rm (law)}}\longrightarrow \int_0^1 B_{2m(t)^2} {\rm d}t,$$
and the distribution of the limit is $\mathcal N(0,4\int_0^1 (1-t)m(t)^2{\rm d}t)$. Moreover, using the fact that $\hat m-m$ is a martingale, we easily prove with Doob's inequality that $\sup_{t\le 1} |\hat m(t)-m(t)|$ converges in probability to 0. Putting all pieces together and applying Slutsky's Lemma gives the result. 

\noindent {\it (ii)} Under ${\bf H_1}$, there exists $t_0\in [0,1]$ such that $\sigma^2(t_0)>m(t_0)$. Then, 
\begin{eqnarray*}
\sqrt{n} \, \hat S_1 & \ge &  \sqrt{n}\big(\hat \sigma^2(t_0)-\hat m(t_0)\big)\\
& \ge & \sqrt{n}\big(\hat \sigma^2(t_0)-\sigma^2(t_0)\big)+\sqrt{n}\big(m(t_0)-\hat m(t_0)\big)+\sqrt{n}\big(\sigma^2(t_0)-m(t_0)\big).
\end{eqnarray*}
The latter term tends to $+\infty$, while the central limit theorem (that can be used because $\mathbb E N_1^4<\infty$) shows that the sequences induced by the first two terms on the right-hand side are stochastically bounded, hence the result with $\hat S_1$. 
Regarding $\hat S_2$, we first observe that under ${\bf H_1}$, 
$\int_0^1 \big(\sigma^2(t)-m(t)\big){\rm d}t >0$, because $\sigma$ and $m$ are right-continuous functions, and $\sigma^2\ge m$. Thus, we only need to prove that the sequences $(\sqrt{n} \int_0 (m(t)-\hat m(t)){\rm d}t)_n$ and $(\sqrt{n} \int_0^1 (\hat \sigma^2(t)-\sigma^2(t)){\rm d}t)_n$ are stochastically bounded. Let us focus on the second sequence. We have
\begin{eqnarray}
\sqrt{n} \int_0^1 \big(\hat \sigma^2(t)-\sigma^2(t)\big){\rm d}t & = & \sqrt{n} \Big(\frac{1}{n-1} \sum_{i=1}^n \int_0^1 \big(N_t^{(i)}\big)^2 {\rm d}t-\mathbb E \int_0^1 N_t^2 {\rm d}t\Big)\nonumber \\
& & - \sqrt{n} \int_0^1 \Big(\frac{n}{n-1} \hat m(t)^2-m(t)^2 \Big) {\rm d}t.\label{sig}
\end{eqnarray}
Since $\mathbb E N_1^4<\infty$, the sequence induced by the first term on the right-hand side is stochastically bounded according to the central limit theorem. Regarding the latter term in \eqref{sig}, we observe that 
\begin{eqnarray*}
\sqrt{n} \Big| \int_0^1 \Big(\frac{n}{n-1}\hat m(t)^2-m(t)^2 \Big){\rm d}t \Big| & \le & 2\sqrt{n} \big(\hat m(1)+m(1)\big) \int_0^1 \big|\hat m(t)-m(t)\big| {\rm d}t\\
& & +m(1)^2.
\end{eqnarray*}
By the Cauchy-Schwarz inequality, there exists a constant $C>0$ such that the $L^1$-norm of the leftmost term is bounded by 
\begin{eqnarray*}
C\Big[1+\sqrt{n}\mathbb E^{1/2} \Big(\int_0^1 \big|\hat m(t)-m(t)\big| {\rm d}t\Big)^2 \Big] & \le &  C\Big[1+\Big(\int_0^1 {\rm var}(N_t){\rm d}t\Big)^{1/2}\Big]\\
& \le &  C\big(1+\mathbb E N_1^2\big).
\end{eqnarray*}
Thus, $(\sqrt{n} \int_0^1 (\hat \sigma^2(t)-\sigma^2(t)){\rm d}t)_n$ is stochastically bounded.
$\Box$ 

\section{Proof of Theorem \ref{local}} 
In this section, we assume that ${\bf H_1^{\it n}}$ holds, hence in particular $N$ is a Cox process with intensity $\lambda=\lambda_0+d_n \Delta$ that depends on $n$. 

For simplicity, we let $Z^{(n)}=(Z_t^{(n)})_{t\in [0,1]}$ the centered process defined for all $t\in [0,1]$ by 
$$Z_t^{(n)}=\sqrt{n} \big(\hat \sigma^2(t)-\sigma^2(t)+m(t)-\hat m(t)\big).$$

\subsection{Auxiliary results} 

\begin{lem} \label{I} Let $k=1,\cdots,6$. Then, there exists $C>0$ independent of $n$ such that for all $s,t\in [0,1]$, $$\mathbb E |N_t-N_s|^k\le C |t-s|^k.$$ \end{lem} 

\noindent {\bf Proof.} Without loss of generality, we assume that $s\le t$. Recall that the distribution of $N_t-N_s$ given $\Lambda$ follows a Poisson distribution with parameter $\Lambda(t)-\Lambda(s)=\int_s^t \lambda(u){\rm d}u$, and that the $k$-th moment of a Poisson distribution with parameter $\mu >0$ is bounded by some constant $C>0$ multiplied by $\mu^k$. Thus, using Jensen's Inequality, we get 
\begin{eqnarray*}
\mathbb E \big(N_t-N_s\big)^k & = & \mathbb E \, \mathbb E \big[\big(N_t-N_s\big)^k|\Lambda\big]\\
& \le & C \mathbb E \Big(\int_s^t \lambda(u){\rm d}u \Big)^k\\
& \le & C (t-s)^{k-1} \int_s^t \mathbb E \lambda^k (u) {\rm d}u \\
& \le & 2^k C \big(\sup_{t\in [0,1]} \lambda_0(t)^k+\sup_{t\in [0,1]} \mathbb E |\Delta_t^k|\big) (t-s)^k,
\end{eqnarray*}
hence the lemma. $\Box$

\begin{lem} \label{EI} For all $t\in [0,1]$, 
$((Z_t^{(n)})^2)_n$ is a uniformly integrable sequence. \end{lem} 

\noindent {\bf Proof.} According to the Rosenthal Inequality and Lemma \ref{I}, there exists a constant $C>0$ that does not depend on $n$ such that
$$n^{3/2}\mathbb E |\hat m(t)-m(t)|^3\le C \mbox{ and } 
n^{3/2}\mathbb E |\hat \sigma^2(t)-\sigma^2(t)|^3 \le C.$$
Thus, we deduce that 
$\sup_n \mathbb E |Z_t^{(n)}|^3<\infty,$
which implies that $((Z_t^{(n)})^2)_n$ is uniformly integrable. $\Box$

\begin{lem} \label{tension} The sequence of processes $(Z^{(n)})_n$ is tight. \end{lem} 

\noindent {\bf Proof.} For an integrable real random variable $Z$, we let $\{Z\}=Z-\mathbb E Z$. First observe that we have, for all $t\in [0,1]$ : 
$$Z_t^{(n)}=\sqrt{n} \Big(\frac{1}{n-1}\mathbb E N_t^2+\frac{1}{n-1} \sum_{i=1}^n \big\{ \big(N_t^{(i)}\big)^2\big\}-\{\hat m(t)\}^2+\{\hat m(t)\}\big(1-2m(t)\big)\Big).$$
We shall make use of the classical criterion of tightness (see e.g. Theorem 13.6 in the book by Billingsley, 1999). Clearly, the sequence will be proved to be tight if if we prove that  each of the sequences of processes defined by
$$X^{1,n}=\sqrt{n}\{\hat m\}, \ X^{2,n}=\sqrt{n}\{\hat m\}^2 \mbox{ and } X^{3,n}=\frac{\sqrt{n}}{n-1} \sum_{i=1}^n \big\{ \big(N_t^{(i)}\big)^2\big\}$$ satisfy the inequality 
\begin{equation} 
\label{eq:52} \
\mathbb E \big(X^{k,n}_t-X_s^{k,n}\big)^2\le C\big(F(t)-F(s)\big)^2,\ \forall\,  0\le s\le t\le 1,
\end{equation} 
for a constant $C>0$ and some nondecreasing and continuous function $F$, both independent of $n$. We only prove it for $X^{3,n}$. In the sequel, $C>0$ is a constant, independent of $n$, and whose value may change from line to line. Observe that 
\begin{eqnarray*}
\mathbb E \big(X_t^{3,n}-X_s^{3,n}\big)^2 & = & \frac{n}{n-1} \mathbb E \big(\big\{N_t^2\big\}-\big\{N_s^2\big\}\big)^2\nonumber\\
& \le & C \mathbb E \big(N_t^2-N_s^2\big)^2+C\big(\mathbb E N_t^2-\mathbb E N_s^2\big)^2 \\
& \le & C \big\{\mathbb E (N_t-N_s)^4\big\}^{1/2}+ C  \mathbb E (N_t-N_s)^2, \label{eq:22}
\end{eqnarray*} 
by Cauchy-Schwarz. Then, Lemma \ref{I} gives 
$$
\mathbb E \big(X_t^{3,n}-X_s^{3,n}\big)^2 \le C (t-s)^2.$$
Consequently, \eqref{eq:52} holds for $k=3$, with the continuous and nondecreasing function $F(t)=t$. $\Box$

\begin{lem} \label{gaus} Let $B=(B_t)_{t\in\mathbb R_+}$ be a real and standard Brownian Motion. Then, if $m_0$ is the function defined for all $t\in [0,1]$ by $m_0(t)=\int_0^t \lambda_0(u){\rm d}u$, we have 
$$Z^{(n)}\stackrel{({\rm law})}\longrightarrow \big(B_{2m_0^2(t)}\big)_{t\in [0,1]}.$$
\end{lem} 

\noindent {\bf Proof.} In the sequel, for $p=1$ or 2, we let 
$$\hat m_p(t)=\frac{1}{n} \sum_{i=1}^n (N_t^{(i)})^p \mbox{ and } m_p(t)=\mathbb E N_t^p.$$
Let $k\ge 1$ and  $0\le t_1<\cdots <t_k\le 1$. According to the central limit theorem for triangular arrays (for instance the Lyapounov condition is easily seen to be true according to Lemma \ref{I}), we know that the $2k$-dimensional random vector defined by 
\begin{displaymath} 
\sqrt{n}\left( \begin{array}{c}
\hat m_1(t_j)-m_1(t_j) \\
\hat m_2(t_j)-m_2(t_j) \\
\end{array} \right)_{j=1,\cdots,k}
\end{displaymath} 
converges to a normal distribution. Now apply the $\delta$-method to deduce that the $3k$-dimensional random vector 
\begin{displaymath} 
\sqrt{n}\left( \begin{array}{c}
\hat m_1(t_j)-m_1(t_j) \\
\hat m_2(t_j)-m_2(t_j) \\
\hat m_1(t_j)^2-m_1(t_j)^2 \\
\end{array} \right)_{j=1,\cdots,k}
\end{displaymath} 
also converges to a normal distribution. Thus, 
$$\sqrt{n} \big(\hat m_2(t_j)-\hat m_1(t_j)^2-m_2(t_j)+m_1(t_j)^2-\hat m_1(t_j)+m_1(t_j)\big)_{j=1,\cdots,k},$$
converges to a $k$-dimensional normal distribution with mean $\mu$ and covariance matrix $\Sigma$, as well as $(Z^{(n)}_{t_j})_{j=1,\cdots,k}$. Since for all $j=1,\cdots,k$, $\hat m(t_j)$ and $\hat \sigma^2(t_j)$ are unbiased estimators of $m(t_j)$ and $\sigma^2(t_j)$, $\mathbb E Z^{(n)}_{t_j}=0$. Thus, by Lemma \ref{EI}, $\mu=0$. We now proceed to compute the variance matrix $\Sigma$. Let $i,j=1,\cdots,k$. With the notation $\{Z\}=Z-\mathbb E Z$ for a integrable real random variable $Z$, we easily see that the difference between 
$
\mathbb E Z_{t_i}^{(n)} Z_{t_j}^{(n)}$ 
and
$$n \mathbb E \big(\{\hat m_2(t_i)\}-\{\hat m_1^2(t_i)\}-\{\hat m_1(t_i)\}\big)\big(\{\hat m_2(t_j)\}-\{\hat m_1^2(t_j)\}-\{\hat m_1(t_j)\}\big)$$
tends to 0 as $n\to\infty$. Moreover, the difference between the latter term and 
$$\frac{1}{n}\mathbb E \Big( \sum_{\ell=1}^n  \{(N_{t_i}^{(\ell)})^2-N_{t_i}^{(\ell)}(2m(t_i)+1)\}\Big)\Big(\sum_{\ell=1}^n \{(N_{t_j}^{(\ell)})^2-N_{t_j}^{(\ell)}(2m(t_j)+1)\}\Big),$$
denoted by $A$, vanishes as well. But, by independence of processes $(N^{(\ell)})_{\ell \le n}$, we have
\begin{eqnarray*}
A & = & \mathbb E \big( \{N_{t_i}^2\}-\{N_{t_i}\}(2m(t_j)+1)\big)\big( \{N_{t_j}^2\}-\{N_{t_i}\}(2m(t_j)+1)\big)\\
& = & {\rm cov} \big(N_{t_i}^2-(2m(t_i)+1)N_{t_i},\, N_{t_j}^2-(2m(t_j)+1)N_{t_j}\big).
\end{eqnarray*}
Recall that, under ${\bf H_1^{\it n}}$, $N$ is a Cox process with intensity $\lambda=\lambda_0+d_n \Delta$. Since $(d_n)_n$ tends to 0, easy calculations prove that, as $n\to\infty$, $A$ converges to 
\begin{equation} 
\label{poisson1} 
{\rm cov} \big(P_{t_i}^2-(2m_0(t_i)+1)P_{t_i},\, P_{t_j}^2-(2m_0(t_j)+1)P_{t_j}\big),
\end{equation} 
where $P$ is a Poisson process with intensity $\lambda_0$. Now use the properties of the Poisson process and the fact that, for a random variable $Q$ that follows a Poisson distribution with parameter $\mu>0$, we have 
$$\mathbb E (Q-\mu)^2=\mathbb E (Q-\mu)^3=\mu \mbox{ and } \mathbb E (Q-\mu)^4=\mu+3\mu^2.$$
Moreover, for all $t\in [0,1]$, $P_t^2-(2m_0(t)+1)P_t+m_0(t)^2=\{P_t\}^2-\{P_t\}$. We easily deduce from the above properties and the independence of the increments of a Poisson process that the covariance in \eqref{poisson1} equals $2m_0(t_i\wedge t_j)^2$. 
Thus, by Lemma \ref{EI}, the $(i,j)$ term of matrix $\Sigma$ is given by the previous formula. The sequence of processes $(Z^{(n)})_n$ being tight according to Lemma \ref{tension}, we have proved that $Z^{(n)}$ converges in distribution to a centered gaussian process $Z$ such that if $s,t\in [0,1]$, $\mathbb E Z_s Z_t=2m_0(s\wedge t)^2$. Such a gaussian process can be written $B_{2m_0^2}$ where $B$ is a standard Brownian Motion on the line, hence the result. $\Box$

\subsection{Proof of Theorem \ref{local}} 

Recall that, by assumption, $v(t)\ne 0$ for some $t\in [0,1]$. Observe that according to \eqref{king}, we have under ${\bf H_1^{\it n}}$, for all $t\in [0,1]$ : 
\begin{eqnarray*}
\sigma^2(t)-m(t)& = & {\rm var} \big(\Lambda(t)\big)\\
& = & {\rm var} \Big(\int_0^t \big(\lambda_0(s)+d_n \Delta_s\big) {\rm d}s\Big)\\
& = & d_n^2 v(t).
\end{eqnarray*}
Thus, 
$$\sqrt{n}\, \hat S_1=\sup_{t\le 1} \big(Z_t^{(n)}+\sqrt{n} d_n^2 v(t)\big).$$
Moreover, it is an easy exercise to prove that $\hat m(T)\to m_0(T)$ and $\hat I\to I_0$, both in probability. If $\sqrt{n} d_n^2\to \infty$, we then have by Lemma \ref{gaus}, 
$$\sqrt{n}\frac{\hat S_1}{\hat m(T)}\, \stackrel{\rm prob.}\longrightarrow \, +\infty,$$
and similarly for $\hat S_2$, hence $(i)$. We now assume that $\sqrt{n}d_n^2\to d<\infty$. By Lemma \ref{gaus} and Slutsky's Lemma, the continuity of the underlying functional gives : 
$$\sqrt{n}\frac{\hat S_2}{\hat I} \, \stackrel{\rm (law)}\longrightarrow \, \frac{1}{I_0} \int_0^1 \big(B_{2m_0(t)^2}+dv(t)\big){\rm d}t.$$
Observing now that the distribution of the latter term equals 
$$2\mathcal N(0,1)+\frac{d}{I_0} \int_0^1 v(t){\rm d} t$$
gives the result for $\hat S_2$. Regarding $\hat S_1$,  
$$\sqrt{n}\frac{\hat S_1}{\hat m(T)}\, \stackrel{\rm (law)}\longrightarrow \, \frac{1}{m_0(T)} \sup_{t\le 1} \big(B_{2m_0(t)^2}+dv(t)\big),$$
hence the theorem. $\Box$

\bigskip\bigskip
\centerline {\Large {\bf References}}

\bigskip

\noindent Billingsley, P. (1999). {\it Convergence of Probability Measures, 2nd Ed.}, Wiley, New-York. 

\smallskip

\noindent Asmussen, S. (2003). {\it Applied Probability and Queues, 2nd Ed.}, Springer, New-York.

\smallskip

\noindent Bj\"{o}rk, T. and Grandell, J. (1988). Exponential inequalities for ruin probabilities in the Cox case, {\it Scandinavian Actuarial Journal}, 77-111. 

\smallskip

\noindent Bohning, D. (1994). A Note on a Test for Poisson Overdispersion, {\it Biometrika}, 418-419.

\smallskip

\noindent Borodin, A.N, Salminen, P. (2002). {\it Handbook of Brownian Motion-Facts and Formulae, 2nd Ed.}, Springer, New-York. 

\smallskip

\noindent Brown, L., Gans, N., Mandelbaum, A., Sakov, A., Shen, H., Zeltyn, S., and Zhao, L. (2005). Statistical analysis of a telephone call center: A queueing-science perspective, {\it Journal of the American statistical association, 100}(469), 36-50.

\smallskip

\noindent Carroll, B. and Ostlie, D. (2007). {\it An Introduction to Modern Astrophysics}, 2nd ed. Benjamin
Cummings, Reading, MA.

\smallskip

\noindent Cox, D.R. and Isham, V. (1980). {\it Point Processes}, Chapman and Hall, London. 

\smallskip

\noindent Denuit, M., Mar\'echal, X., Pitrebois, S., and Walhin, J. F. (2007). {\it Actuarial modelling of claim counts: Risk classification, credibility and bonus-malus systems}, John Wiley $\&$ Sons.

\smallskip

\noindent Gerstner, W. and Kistler, W. (2002). {\it Spiking Neuron Models: Single Neurons, Populations,
Plasticity}, Cambridge Univ. Press, Cambridge.

\smallskip

\noindent Grandell, J. (1991). {\it Aspects of Risk Theory}, Springer-Verlag, New-York. 

\smallskip

\noindent Heuer, A., Mueller, C., and Rubner, O. (2010). Soccer: Is scoring goals a predictable Poissonian process ? {\it EPL (Europhysics Letters), 89}(3), 380-7.

\smallskip

\noindent Jacod, J. and Shiryaev, A.N. (2003). {\it Limit Theorems for Stochastic Processes, 2nd Ed.}, Springer, New-York. 

\smallskip

\noindent Karr, A.F. (1991). {\it Point Processes and their Statistical Applications, 2nd Ed.}, Marcel Dekker, New-York. 

\smallskip

\noindent Kingman, J.F.C. (1993). {\it Poisson Processes}, Oxford Studies in Probability, Oxford. 

\smallskip

\noindent Kou, S. C. (2008). Stochastic networks in nanoscale biophysics: Modeling enzymatic reaction of
a single protein, {\it J. Amer. Statist. Assoc.}, 961-5.

\smallskip

\noindent Kou, S. C., Xie, X. S. and Liu, J. S. (2005). Bayesian analysis of single-molecule experimental
data (with discussion), {\it J. Roy. Statist. Soc. Ser. C}, 469-6.

\smallskip

\noindent Mandelbaum, A., Sakov, A., and Zeltyn, S. (2000). Empirical analysis of a call center, {\it URL http://iew3. technion. ac. il/serveng/References/ccdata. pdf. Technical Report.}

\smallskip

\noindent Rao, C. R., and Chakravarti, I. M. (1956). Some small sample tests of significance for a Poisson distribution, {\it Biometrics, 12}(3), 264-282.

\smallskip

\noindent Reynaud-Bourret, P., Rivoirard, V., Grammont , F., Tuleau-Malot, C. (2014). Goodness-of-fit tests and nonparametric adaptive estimation for spike train analysis, {\it Journal of Mathematical Neuroscience}, 4:3.

\smallskip

\noindent Revuz, D. and Yor, M. (1999). {\it Continuous Martingales and Brownian Motion, 3rd Ed.}, Springer, New-York. 

\smallskip

\noindent Scargle, J. D. (1998). Studies in astronomical time series analysis. V. Bayesian blocks, a new
method to analyze structure in photon counting data, {\it Astrophys. J.}, 405-8.

\smallskip

\noindent Schmidili, H. (1996). Lundberg inequalities for a Cox model with a piecewise constant intensity, {\it J. Applied Probability}, 196-210. 

\smallskip

\noindent van der Vaart, A.W. (1998). {\it Asymptotic Statistics}, Cambridge Series in Statistical and Probabilistic Mathematics. 

\smallskip


\noindent Zhang, T. and Kou, S.C. (2010). Nonparametric inference of doubly stochastic Poisson process data via the kernel method, {\it Ann. Applied Statist.}, 1913-1941. 

\end{document}